\documentclass{amsart}
\usepackage[foot]{amsaddr}

\usepackage[linktocpage=true]{hyperref}
\hypersetup{
    hypertexnames=true,
    colorlinks=true,
    linkcolor=blue,  
    urlcolor=blue,
    citecolor=blue
}

\usepackage{amssymb}
\usepackage[dvipsnames,table,xcdraw]{xcolor}
\usepackage{amsmath,amsthm}
\usepackage{enumerate}
\usepackage{tikz-cd,graphicx,float}
\usetikzlibrary{arrows}
\usepackage{tikz}
\usetikzlibrary{cd,arrows,shapes.misc}
	\tikzset{thick/.style={line width=.4mm}}
	\tikzstyle{dot}=[circle,thick,fill=red!60]
	\tikzstyle{tinydot}=[circle,thick,fill=red!60,inner sep=.7mm]

\newcommand{\End}{\operatorname{End}}
\DeclareMathOperator{\id}{id}

\newtheorem{Theorem}{Theorem}
\newtheorem*{Theorem*}{Theorem}

\newtheorem*{Corollary*}{Corollary}
\newtheorem{Proposition}{Proposition}

\theoremstyle{definition}
\newtheorem{Definition}{Definition}
\theoremstyle{definition}
\newtheorem{Remark}{Remark}
\theoremstyle{remark}
\newtheorem{Example}{Example}

\let\emph\relax
\DeclareTextFontCommand{\emph}{\bfseries}

\newcommand{\cat}[1]{{\normalfont\textsf{#1}}}

\def\O{\mathcal{O}}
\def\C{\cat{C}}

\def\c{\colon}
\def\o{\circ}

\def\R{\mathbb{R}}

\title{Entropy as a topological operad derivation}
\author{Tai-Danae Bradley}
\address{Sandbox@Alphabet, Mountain View, CA 94043}
\email{tai.danae@math3ma.com}

\begin{document}
\begin{abstract}
We share a small connection between information theory, algebra, and topology---namely, a correspondence between Shannon entropy and derivations of the operad of topological simplices. We begin with a brief review of operads and their representations with topological simplices and the real line as the main example. We then give a general definition for a derivation of an operad in any category with values in an abelian {bi}module over the operad. The main result is that Shannon entropy defines a derivation of the operad of topological simplices, and that for every derivation of this operad there exists a point at which it is given by a constant multiple of Shannon entropy. We show this is compatible with, and relies heavily on, a well-known characterization of entropy given by Faddeev in 1956 and a recent variation given by Leinster.
\end{abstract}
\maketitle

\section{Introduction}\label{sec:intro}
In this article, we describe a simple connection between information theory, algebra, and topology. To motivate the idea, consider the function $d\colon[0,1]\to\mathbb{R}$ defined by

\begin{equation*}
d(x) = 
\begin{cases}
-x\log x & \text{if $x>0$},\\
0 & \text{if $x=0$}.
\end{cases}
\end{equation*}

This map satisfies an equation reminiscent of the Leibniz rule from Calculus, $d(xy)=d(x)y+xd(y)$ for all $x,y\in[0,1]$. In other words, $d$ is a nonlinear derivation \cite{leinster2021entropy}, (Lemma 2.2.6). This derivation may also bring to mind the Shannon entropy of a probability distribution. Indeed, a probability distribution on a finite set $\{1,\ldots,n\}$ for $n\geq 1$ is a tuple {of nonnegative real numbers} $p=(p_1,\ldots,p_n)$ satisfying $\sum_{i=1}^n p_i=1$, and the \emph{Shannon entropy} of $p$ is defined to be

\[H(p)= -\sum_{i=1}^np_i\log p_i=\sum_{i=1}^n d(p_i).\]
Although $d$ is not linear, this may prompt one to wonder about settings in which Shannon entropy itself is a derivation. We describe one such setting below by showing a correspondence between Shannon entropy and derivations of the operad of topological simplices. 

\subsection{Motivation}
As evidenced by recent work, the intersection of information theory and algebraic topology is fertile ground. In 2015 tools of information cohomology were introduced \mbox{in \cite{baudot2015}} by Baudot and Bennequin who construct a certain cochain complex for which entropy represents the unique cocycle in degree 1. In the same year, Elbaz-Vincent and Gangl approached entropy from an algebraic perspective and showed that {what are known as} information functions of degree 1 behave ``a lot like certain derivations'' \cite{poly2015}. A few years prior in 2011, Baez, Fritz, and Leinster gave a category theoretical characterization of entropy in \cite{BFL}, which was recently extended to the quantum setting by Parzygnat in \cite{arthur}. In preparation of that 2011 result, Baez remarked in the informal article \cite{EFunctor} that entropy appears to behave similarly to a derivation in a certain operadic context, an observation we verify and make explicit below. Cohomological ideas are also explored in Mainiero's recent work, where entropy is found to appear in the Euler characteristic of a particular cochain complex associated to a quantum state \cite{Mainiero}. Upon taking inventory, one thus has the sense that entropy behaves somewhat similar to ``$d$ of something,'' for some (co)boundary-like operator $d.$ The present article is in this same vein. Notably, once a few simple definitions are in place, the mathematics is quite straightforward. Even so, we feel it is worth sharing if for no other reason than to provide a glimpse at yet another algebraic and topological facet of entropy.

\subsection{Background} To start, our work is based on a particular characterization of Shannon entropy that is compatible with an operadic viewpoint. Let $\Delta^n$ denote the standard {topological $n$-simplex} for $n\geq 0$, 

\[\Delta^n:=\{(p_0,p_1,\ldots,p_n)\in\mathbb{R}^{n+1}\mid 0\leq p_i\leq 1 \text{ and }\sum_{i=0}^np_i=1\},\]
where $\Delta^0$ denotes the unique probability distribution on the one-point set. More generally, any probability distribution $p=(p_0,\ldots,p_n)$ on an $n+1$-element set is a point in $\Delta^n$. Given $n+1$ probability distributions $q^i=(q^i_0,\ldots, q^i_{k_i})\in\Delta^{k_i}$ where $i=0,1,\ldots,n$, they may be composed with $p$ simultaneously to obtain a point in {$\Delta^{k_0+k_1+\cdots+k_n+n}$} denoted by
\[p\circ (q^0,q^1,\ldots,q^n):=(p_0q^0_0,\ldots,p_0q^0_{k_0},p_1q^1_1,\ldots,p_1q^1_{k_1},\ldots,p_nq^n_{1},\ldots,p_nq^n_{k_n}).\]
As shown in \cite{leinster2021entropy} and reviewed below, this composition of probabilities finds a natural home in the language of operads. Furthermore,  it plays a key role in a well-known 1956 characterization of Shannon entropy due to {D. K.} Faddeev \cite{Fad}. A proof of a slight variation of Faddeev's result was recently given by Leinster \cite{leinster2021entropy}, (Theorem 2.5.1). That is the version we quote here.

\begin{Theorem}[Faddeev-Leinster]\label{thm:FL}
Let $\{F\colon \Delta^n\to\mathbb{R}\}_{n\geq 0}$ be a sequence of functions. The following are equivalent:
\begin{enumerate}
	\item the functions F are continuous and satisfy
	
	\begin{equation}\label{eq:Faddeev}
	F(p\circ (q^0,\ldots,q^n)) = F(p) + \sum_{i=0}^np_iF(q^i)
	\end{equation}
	
	where $n\geq 0$ and $p\in\Delta^n$ and $q^i\in\Delta^{k_i}$ with {$k_0$}, $k_1,\ldots,k_n\geq 0$;
	\item $F=cH$ for some $c\in\mathbb{R}.$
\end{enumerate}
\end{Theorem}

\noindent To make the connection with derivations, let us introduce some notation. Given a probability distribution $p\in\Delta^n$ let $\bar{p} \colon \mathbb{R}^{n+1}\to \mathbb{R}$ denote the function that maps a point $x=(x_0,\ldots,x_n)$ to the standard inner product $\langle p,x\rangle=\sum_{i=0}^np_ix_i$. Then, when $F=H$, Equation (\ref{eq:Faddeev}) may be rewritten as

\begin{equation}\label{eq:Baez}
H(p\circ (q^0,\ldots,q^n)) = H(p) + \bar p(H(q^0),\ldots , H(q^n)).
\end{equation}
This equation is one hint that entropy might be a derivation, although a ``$q$'' is notably absent from the first term on the right-hand side. As a further teaser, Baez explored an algebraic interpretation of Equation (\ref{eq:Baez}) in the {informal} article \cite{EFunctor}, where the reader is reminded that Shannon entropy is a derivative of the partition function of a probability distribution with respect to Boltzmann's constant, considered as a formal parameter. In that article, Equation (\ref{eq:Baez}) follows in a few short lines from this computation. One is thus motivated to look for a general framework of operad derivations for which Equation (\ref{eq:Baez}) is an example. This is what we describe below.

Section \ref{sec:definitions} reviews the definition of operads and representations of them. We will recall that {the collection of topological simplices admits the structure of an operad} as in \cite{leinster2021entropy} and that $\mathbb{R}$ gives rise to a representation of it. In Section \ref{sec:mainresult}, we define an abelian {bi}module $M$ over any operad $\mathcal{O}$ and the notion of a derivation of $\mathcal{O}$ with values in $M$. With these definitions in place, Equation (\ref{eq:Baez}) will find a generalization in Proposition \ref{proposition}, and the main result will quickly follow{.}

\begin{Theorem*} Shannon entropy defines a derivation of the operad of topological simplices, and for every derivation of this operad there exists a point at which it is given by a constant multiple of Shannon entropy.
\end{Theorem*}

\section{Acknowledgements}
I thank Darij Grinberg, Joey Hirsh, Tom Leinster, Jim Stasheff, and John Terilla for helpful discussions as well as the anonymous referees for their insightful feedback.

\section{Background: Operads and Their Representations}\label{sec:definitions}

\begin{figure}[H]
\begin{equation*}
\begin{tikzpicture}[x=.5cm,y=.5cm,baseline={(current bounding box.center)},scale=0.9]
	\node[] (a) at (0,0) {};
	\node[] (l) at (1,1) {};
	\node[] (r) at (-1,1) {};
	\node[] (m) at (0,1) {};
	\draw[thick, shorten <=-2mm, shorten >=-1mm, color=magenta] (0,-.75) -- (a);
	\draw[thick, shorten <=-2mm, color=magenta] (a) -- (l);
	\draw[thick, shorten <=-2mm, color=magenta] (a) -- (r);
	\draw[thick, shorten <=-2mm, color=magenta] (a) -- (m);
	\node [] at (0,-2) {$f$};
\end{tikzpicture}
\qquad
\begin{tikzpicture}[x=.6cm,y=.65cm,baseline={(current bounding box.center)},scale=0.9]
	\node[] (a) at (0,0) {};
	\node[] (1) at (-.86,1/2) {};
	\node[] (2) at (-1/2,.86) {};
	\node[] (3) at (1/2,.86) {};
	\node[] (4) at (.86,1/2) {};
	\draw[thick, shorten <=-1.5mm, shorten >=-1.5mm, color=Cerulean] (0,-.75) -- (a);
	\draw[thick, shorten <=-1.5mm, color=Cerulean] (a) -- (1);
	\draw[thick, shorten <=-1.5mm, color=Cerulean] (a) -- (2);
	\draw[thick, shorten <=-1.5mm, color=Cerulean] (a) -- (3);
	\draw[thick, shorten <=-1.5mm, color=Cerulean] (a) -- (4);
	\node [] at (0,-1.5) {$g$};
\end{tikzpicture}
\quad\rightsquigarrow\quad
\begin{tikzpicture}[x=.5cm,y=.5cm,baseline={(current bounding box.center)},scale=0.9]
	\node[] (a) at (0,0) {};
	\node[] (l) at (1,1) {};
	\node[] (r) at (-1,1) {};
	\node[] (m) at (0,1) {};
	\draw[thick, shorten <=-1.5mm, shorten >=-1mm, color=magenta] (0,-.75) -- (a);
	\draw[thick, shorten <=-1.5mm, color=magenta] (a) -- (l);
	\draw[thick, shorten <=-1.5mm, color=magenta] (a) -- (r);
	\draw[thick, shorten <=-1.5mm, color=magenta] (a) -- (m);

	\node[] (T) at (0,1.5) {};
	\draw[thick,color=Cerulean,shorten <=-1.5mm,shorten >=-1mm] (0,1) -- (T);
	\draw[thick,color=Cerulean,shorten <=-1.5mm] (T) -- (-.86,2);
	\draw[thick,color=Cerulean,shorten <=-1.5mm] (T) -- (-1/2,2.36);
	\draw[thick,color=Cerulean,shorten <=-1.5mm] (T) -- (1/2,2.36);
	\draw[thick,color=Cerulean,shorten <=-1.5mm] (T) -- (.86,2);

	\node [] at (0,-2) {$f\circ_2g$};
\end{tikzpicture}
\end{equation*}
\caption{One of the three ways to compose a $4$-ary operation $g$ with a $3$-ary operation $f$.}
\label{fig:operads}
\end{figure}
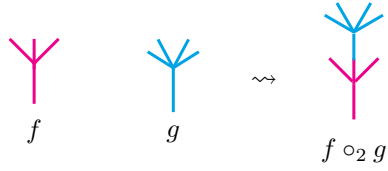

In an introduction to operads, it is helpful to first think about algebras. An algebra $A$ is a vector space $V$ equipped with a bilinear map $\mu\colon V\times V\to V$ thought of as multiplication. Depending on whether $\mu$ satisfies a particular {relation}, the algebra {will usually be described by an approriate qualifier}. For instance, if $\mu(v,w)=\mu(w,v)$ for all $v,w\in V$, then $A$ is called a \textit{commutative algebra}; if $\mu(\mu(u,v),w)=\mu(u,\mu(v,w))$ for all $u,v,w\in V$, then $A$ is a called an \textit{associative algebra}, and so on. Behind each of these algebras is a particular operad that encodes the behavior of the multiplication map $\mu$. To motivate the formal definition, it is helpful to visualize $\mu$ as a planar binary rooted tree and more generally to imagine an arbitrary $n$-ary operation as a planar rooted tree with $n$ leaves. There is a natural way to compose such operations. For instance, when $f$ is a $3$-ary operation and $g$ is a $4$-ary operation, they may be composed to obtain a $6$-ary operation by using the output of $g$ as one of the inputs of $f$ as illustrated in Figure~\ref{fig:operads}. There $g$ has been grafted into the \textit{second} leaf of the tree associated to $f$, and so we denote that choice with the subscript ``$\circ_2$'' in the figure. There are two other composites $f\circ_1g$ and $f \circ_3 g$, which are not shown but are obtained similarly.

In general, there are $n$ ways to compose an $m$-ary operation with an $n$-ary operation, and the resulting operation will always have arity $m+n-1$. This composition should further satisfy some sensible associativity and unital axioms, and the collection of all such operations with their compositions is called an \emph{operad}. The concept has origins in category theory \cite{lambek} and has been used extensively in algebraic topology and homotopy \mbox{theory \cite{may2006geometry,BV73,loday2012algebraic,Bruno,Stasheff}} with applications in physics as well \cite{Mar,markl2002operads}. Operads may be defined in any {symmetric monoidal} category, and for ease of exposition below, we will assume all categories $\mathsf{C}$ are concrete {(that is, all objects have underlying sets)} so that we may refer to \textit{elements} in a given object of $\mathsf{C}$. Indeed, the main example to have in mind is the category of topological spaces.

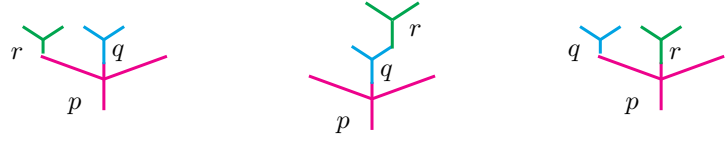
\begin{figure}[H]
\begin{equation*}
\begin{tikzpicture}[x=.55cm,y=.5cm,baseline={(current bounding box.center)},scale=0.7]
	\node[] (a) at (0,0) {};
	\node[] (r) at (2.5,1) {};
	\node[] (l) at (-2.5,1) {};
	\node[] (m) at (0,1) {};
	\draw[thick, shorten <=-1.5mm, shorten >=-1mm, color=magenta] (0,-.75) -- (a);
	\draw[thick, shorten <=-1.5mm, color=magenta] (a) -- (r);
	\draw[thick, shorten <=-1.5mm, color=magenta] (a) -- (l);
	\draw[thick, shorten <=-1.5mm, color=magenta] (a) -- (m);
	
	\node[] (T) at (0,1.5) {};
	\draw[thick,color=Cerulean,shorten <=-1.5mm,shorten >=-1.2mm] (0,1) -- (T);
	\draw[thick,color=Cerulean,shorten <=-1.5mm] (T) -- (-.7,2);
	\draw[thick,color=Cerulean,shorten <=-1.5mm] (T) -- (.7,2);

	\node[] (p1) at (-2.1,1.5) {};
	\draw[thick,color=Green,shorten <=-1.1mm,shorten >=-.6mm] (-2.1,1.2) -- (p1);
	\draw[thick,color=Green,shorten <=-1.5mm] (p1) -- (-2.8,2);
	\draw[thick,color=Green,shorten <=-1.5mm] (p1) -- (-1.4,2);

	\node[] at (-1,-1) {$p$};
	\node[] at (-3,1) {$r$};
	\node[] at (.5,1) {$q$};
\end{tikzpicture}
\qquad\qquad
\begin{tikzpicture}[x=.55cm,y=.5cm,baseline={(current bounding box.center)},scale=0.7]
	\node[] (a) at (0,0) {};
	\node[] (r) at (2.5,1) {};
	\node[] (l) at (-2.5,1) {};
	\node[] (m) at (0,1) {};
	\draw[thick, shorten <=-1.5mm, shorten >=-1mm, color=magenta] (0,-.75) -- (a);
	\draw[thick, shorten <=-1.5mm, color=magenta] (a) -- (r);
	\draw[thick, shorten <=-1.5mm, color=magenta] (a) -- (l);
	\draw[thick, shorten <=-1.5mm, color=magenta] (a) -- (m);
	
	\node[] (T) at (0,1.5) {};
	\draw[thick,color=Cerulean,shorten <=-1.5mm,shorten >=-1.2mm] (0,1) -- (T);
	\draw[thick,color=Cerulean,shorten <=-1.5mm] (T) -- (-.7,2);
	\draw[thick,color=Cerulean,shorten <=-1.5mm] (T) -- (.7,2);

	\node[] (p) at (.7,3) {};
	\draw[thick,color=Green,shorten <=-1.1mm,shorten >=-.3mm] (p) -- (.7,2);
	\draw[thick,color=Green,shorten <=-1.5mm,shorten >=-1mm] (p) -- (1.4,3.5);
	\draw[thick,color=Green,shorten <=-1.5mm,shorten >=-1mm] (p) -- (0,3.5);

	\node[] at (-1,-1) {$p$};
	\node[] at (.5,1) {$q$};
	\node[] at (1.5,2.6) {$r$};
\end{tikzpicture}
\qquad\qquad
\begin{tikzpicture}[x=.55cm,y=.5cm,baseline={(current bounding box.center)},scale=0.7]
	\node[] (a) at (0,0) {};
	\node[] (r) at (2.5,1) {};
	\node[] (l) at (-2.5,1) {};
	\node[] (m) at (0,1) {};
	\draw[thick, shorten <=-1.5mm, shorten >=-1mm, color=magenta] (0,-.75) -- (a);
	\draw[thick, shorten <=-1.5mm, color=magenta] (a) -- (r);
	\draw[thick, shorten <=-1.5mm, color=magenta] (a) -- (l);
	\draw[thick, shorten <=-1.5mm, color=magenta] (a) -- (m);
	
	\node[] (T) at (0,1.5) {};
	\draw[thick,color=Green,shorten <=-1.5mm,shorten >=-1.2mm] (0,1) -- (T);
	\draw[thick,color=Green,shorten <=-1.5mm] (T) -- (-.7,2);
	\draw[thick,color=Green,shorten <=-1.5mm] (T) -- (.7,2);

	\node[] (p1) at (-2.1,1.5) {};
	\draw[thick,color=Cerulean,shorten <=-1.1mm,shorten >=-.6mm] (-2.1,1.2) -- (p1);
	\draw[thick,color=Cerulean,shorten <=-1.5mm] (p1) -- (-2.8,2);
	\draw[thick,color=Cerulean,shorten <=-1.5mm] (p1) -- (-1.4,2);

	\node[] at (-1,-1) {$p$};
	\node[] at (-3,1) {$q$};
	\node[] at (.5,1) {$r$};
\end{tikzpicture}
\end{equation*}
\caption{Associativity in an operad. Left) First composing $q$ with $p$ and then $r$ is the same as first composing $r$ with $p$ and then $q$. The order in which this is performed does not matter. (Right) The same is true if $r$ appears to the right, rather than the left, of $q.$ (Middle) Likewise, $r$ may first be composed with $q$ and their composite may then be composed with $p$, or $q$ may be first composed with $p$ followed by $r.$ Again, the order does not matter.}
\label{fig:operad_axioms}
\end{figure}

\begin{Definition} Let $\mathsf{C}$ be a {symmetric monoidal} category {with monoidal product $\otimes.$}. An \emph{operad} in $\mathsf{C}$ consists of a sequence of objects $\{\O(1), \O(2),\ldots\}$ together with morphisms

	\[\o_i\c \O(n)\otimes \O(m)\to \O(n+m-1)\]	
in $\mathsf{C}$ for all $n,m\geq 1$ and $1\leq i \leq n$ and an operation $1\in \O(1)$ satisfying the following:
\begin{itemize}
	\item[(i)][associativity] For all $p\in\O(n)$ and $q\in \O(m)$ and $r\in \O(k)$, 
	
\begin{align*}
 (p\o_j q)\o_i r =
 \begin{cases}
 	(p\o_i r)\o_{j+k-1} q &\text{if $1\leq i \leq j-1$}\\
 	p\o_j(q\o_{i-j+1}r) &\text{if $j\leq i \leq j+m-1$}\\
 	(p\o_{i-m+1}r)\o_j q &\text{if $i \geq j+m$}
 \end{cases}
 \end{align*}
 
	\item[(ii)][identity] The operation $1\in \O(1)$ acts as an identity in the sense that
	
 \[1\circ_1 p = p\circ_i 1=p\]
 
 for all $p\in \O(n)$ and $1\leq i \leq n.$
 \end{itemize}
\end{Definition}
The definition is conceptually simple despite its cumbersome appearance. For instance, Figure~\ref{fig:operad_axioms} illustrates the associativity requirements listed in item (i). 

As mentioned above, one often thinks of the elements $\O(n)$ as abstract $n$-to-1 operations, and the morphisms $\circ_i$ specify a way to compose them. {It is common to begin indexing the sequence of objects at $n=0$ to account for $0$-ary operations, but as we will soon see, our main example of an operad in Example} \ref{ex:simplices} {will have no $0$-ary operations, and so our definition starts with $\O(1)$.} We do not consider an action of the symmetric group and so $\O$ is sometimes called a \textit{non-symmetric operad}, but we will simply call it an operad. In the special case when $\C$ is the category of vector spaces with linear maps {and $\otimes$ is the tensor product}, $\O$ is often called a \textit{linear operad}. When it is the category $\mathsf{Top}$ of topological spaces with continuous maps {and $\otimes$ is the Cartesian product}, $\O$ is often called a \textit{topological operad}.

\begin{Example}\label{ex:endomorphism}
{Given a set $X$,} the \emph{endomorphism operad} is $\End_X=\{\End_X(1),\End_X(2),\ldots\}$ where $\End_X(n):=\C(X^n,X)$ denotes the set of all {functions} from {the $n$-fold Cartesian product} $X^n$ to $X$. The unit operation in $\End_X(1)$ is the identity {function} $\id_X\colon X\to X.$ If $f\in \C(X^n,X)$ and $g\in \C(X^m,X)$ are a pair of {functions}, then for each $i=1,\ldots,n$ the composition $f\circ_i g$ is obtained by using the output of $g$ as the $i$th input of $f.$ Explicitly, given $(x_1,\ldots,x_{n+m-1})\in~ X^{n+m-1}$,

\[(f\circ_i g)(x_1,\ldots,x_{n+m-1}):=f(x_1,\ldots,x_{i-1},g(x_i,\ldots,x_{i+m-1}),x_{i+m},\ldots, x_{n+m-1}).\]

{The simultaneous composition of several functions may also be considered. That is, given $n$ functions $g_i\in\mathsf{C}(X^{k_i},X)$ where $i=1,\ldots,n$ they may be composed with $f$ simultaneously to obtain a new function $f\circ (g_1,\ldots,g_n)\in\mathsf{C}(X^{k_1+\cdots+k_n},X)$, which is again defined by using the outputs of the $g_i$ as the inputs of $f.$ Explicitly, given $(x_1,\ldots,x_{k_1+\cdots+k_n})\in X^{k_1+\cdots+k_n}$, we have}
\begin{equation*}
    (f\circ (g_1,\ldots,g_n))(x_1,\ldots,x_{k_1+\cdots+k_n})=f(g_1(x_1,\ldots,x_{k_1}),\ldots,g_n(x_{k_1+\cdots+k_{n-1}+1},\ldots,x_{k_1+\cdots+k_n}))
\end{equation*}
\end{Example}

\begin{Example}\label{ex:simplices} The simplices $\Delta^0,\Delta^1,\Delta^2,\ldots$ give rise to a topological operad  called \emph{the operad of topological simplices}  $\Delta=\{\Delta_1,\Delta_2,\ldots\}$ where $\Delta_n:=\Delta^{n-1}$. The unit operation in $\Delta_1$ is the unique probability distribution on a one-point set. If $p=(p_1,\ldots,p_n)\in \Delta_n$ and $q=(q_1,\ldots,q_m)\in \Delta_m$ are probability distributions, then the composition $p\circ_i q$ is obtained by multiplying each of the $m$ coordinates of $q$ by $p_i$ and then replacing the $i^{\text{th}}$ coordinate of $p$ with the resulting $m$-tuple. Explicitly,

\[p\circ_i q:=(p_1,\ldots,p_iq_1,\ldots,p_iq_m,\ldots,p_n)\in\Delta_{n+m-1}.\]
Equivalently, the distribution $p$ may be visualized as a planar tree with $n$ leaves labeled by the probabilities $p_1,\ldots,p_n$ and similarly for $q$. Then the composition $p\circ_i q$ is obtained by ``painting'' each of the leaves of $q$ with the probability $p_i$ and grafting the resulting tree into the $i^{\text{th}}$ leaf of $p$ as below. Notice the sum of {the probabilities on} the leaves on the composite tree is 1.

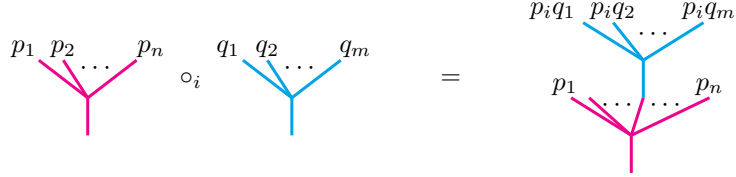
\[
	\begin{tikzpicture}[x=.65cm,y=.5cm,baseline={(current bounding box.center)}]
	\draw[thick,color=magenta] (0,0) --(-1,1); \node at (-1.3,1.3) {$p_1$};
	\draw[thick,color=magenta] (0,0) --(-.5,1); \node at (-.5,1.3) {$p_2$};
	\node at (0.2,0.8) {$\cdots$};
	\draw[thick,color=magenta] (0,0) --(1,1); \node at (1.3,1.3) {$p_n$};
	\draw[thick,color=magenta] (0,-1) --(0,0);
	\end{tikzpicture}
	\circ_i
		\begin{tikzpicture}[x=.65cm,y=.5cm,baseline={(current bounding box.center)}]
	\draw[thick,color=Cerulean] (0,0) --(-1,1); \node at (-1.3,1.3) {$q_1$};
	\draw[thick,color=Cerulean] (0,0) --(-.5,1); \node at (-.5,1.3) {$q_2$};
	\node at (0.2,0.8) {$\cdots$};
	\draw[thick,color=Cerulean] (0,0) --(1,1); \node at (1.3,1.3) {$q_m$};
	\draw[thick,color=Cerulean] (0,-1) --(0,0);
	\end{tikzpicture}
\qquad = \qquad
	\begin{tikzpicture}[x=.8cm,y=.5cm,baseline={(current bounding box.center)}]
	\draw[thick,color=magenta] (0,0) --(-1,1); \node at (-1.1,1.3) {$p_1$};
	\draw[thick,color=magenta] (0,0) --(-.7,1); 
	\node at (-.2,0.8) {$\cdots$};
	\draw[thick,color=magenta] (0,0) --(.2,1); 
	\node at (.6,0.8) {$\cdots$};
	\draw[thick,color=Cerulean] (.2,1)--(.2,2);
	\draw[thick,color=Cerulean] (.2,2)--(-.8,3); \node at (-1.3,3.3) {$p_iq_1$};
	\draw[thick,color=Cerulean] (.2,2)--(-0.3,3);\node at (-0.3,3.3) {$p_iq_2$};
	\node at (.4,2.7) {$\cdots$};
	\draw[thick,color=Cerulean] (.2,2)--(1.2,3);\node at (1.3,3.3) {$p_iq_m$};
	\draw[thick,color=magenta] (0,0) --(1.3,1); \node at (1.3,1.3) {$p_n$};
	\draw[thick,color=magenta] (0,-1) --(0,0);
	\end{tikzpicture}
\]
As an example, if $p=\left(\frac{1}{6},\ldots,\frac{1}{6}\right)$ represents the probability distribution of rolling a six-sided die and $q=\left(\frac{1}{2},\frac{1}{2}\right)$ is that of a fair coin toss, then $p\circ_3 q=\left(\frac{1}{6}, \frac{1}{6},\frac{1}{12},\frac{1}{12}, \frac{1}{6},\frac{1}{6},\frac{1}{6}\right)$ is a point in $\Delta_7$, whose picture is shown on the left of Figure~\ref{fig:dice_pic}.

\begin{figure}[H]
\begin{equation*}
	\begin{tikzpicture}[x=.65cm,y=.65cm,baseline={(current bounding box.center)}]
	\draw[thick,color=magenta] (0,0) -- (-3,1);
		\node at(-3.1,1.4) {$\frac{1}{6}$}; 
	\draw[thick,color=magenta] (0,0) -- (-1.8,1); 
		\node at (-1.8,1.4){$\frac{1}{6}$}; 
	\draw[thick,color=magenta] (0,0) -- (-0.6,1);
		\draw[thick,color=Cerulean] (-0.6,1) -- (-0.6,1.8);
		\draw[thick,color=Cerulean] (-0.6,1.8) -- (-1.2,2.4);
			\node at (-1.2,2.8){$\frac{1}{12}$};
		\draw[thick,color=Cerulean] (-0.6,1.8) -- (0,2.4);
			\node at (0,2.8){$\frac{1}{12}$};
	\draw[thick,color=magenta] (0,0) -- (0.6,1);
		\node at (0.7,1.4) {$\frac{1}{6}$};
	\draw[thick,color=magenta] (0,0) -- (1.8,1);
		\node at (1.9,1.4){$\frac{1}{6}$};
	\draw[thick,color=magenta] (0,0) --(3,1); \node at (3.1,1.4){$\frac{1}{6}$};
	\draw[thick,color=magenta] (0,-1) --(0,0);
	\end{tikzpicture}
\qquad\qquad
	\begin{tikzpicture}[x=1cm,y=.5cm,baseline={(current bounding box.center)}]
	\draw[thick,color=magenta] (0,0) -- (-3,1); 
		\draw[thick,color=Green] (-3,1)--(-3,2); 
		\draw[thick,color=Green] (-3,2)--(-3.5,2.5); \node at (-3.8,2.8){$p_1q^1_1$};
		\draw[thick,color=Green] (-3,2)--(-3.3,2.5);
		\draw[thick,color=Green] (-3,2)--(-3.1,2.5);
		\node at (-2.9,2.3) {$\cdot\cdot$};
		\draw[thick,color=Green] (-3,2)--(-2.5,2.5); \node at (-2.4,2.8){$p_1q^1_{k_1}$};
	\draw[thick,color=magenta] (0,0) -- (-.3,1);  
		\draw[thick,color=Cerulean] (-.3,1)--(-.3,2); 
		\draw[thick,color=Cerulean] (-.3,2)--(.2,2.5); \node at (-1.1,2.8){$p_2q^2_1$ };
		\draw[thick,color=Cerulean] (-.3,2)--(-.7,2.5);
		\draw[thick,color=Cerulean] (-.3,2)--(-.5,2.5);
		\node at (-.2,2.3) {$\cdot\cdot$};
		\draw[thick,color=Cerulean] (-.3,2)--(-.9,2.5);\node at (.3,2.8){$p_2q^2_{k_2}$};
	\node at (.8,0.8) {$\cdots$}; 
	\draw[thick,color=magenta] (0,0) --(2.5,1); 
		\draw[thick,color=Dandelion] (2.5,1)-- (2.5,2); 
		\draw[thick,color=Dandelion] (2.5,2)--(2,2.5);\node at (1.7,2.8){$p_nq^n_1$};
		\draw[thick,color=Dandelion] (2.5,2)--(2.2,2.5);
		\draw[thick,color=Dandelion] (2.5,2)--(2.4,2.5);
		\node at (2.6,2.3) {$\cdot\cdot$};
		\draw[thick,color=Dandelion] (2.5,2)--(3,2.5);\node at (3.3,2.8){$p_nq^n_{k_n}$};
	\draw[thick,color=magenta] (0,-1) --(0,0);
	\end{tikzpicture}
\end{equation*}
\caption{(Left) A picture of the composition $p\circ_3q$ when $p$ is the probability distribution associated to a six-sided die and $q$ is that of a fair coin toss. (Right) The simultaneous composition of $n$ probability distributions $q^i\in\Delta_{k_i}$ with a given $p\in\Delta_n$.}
\label{fig:dice_pic}
\end{figure}
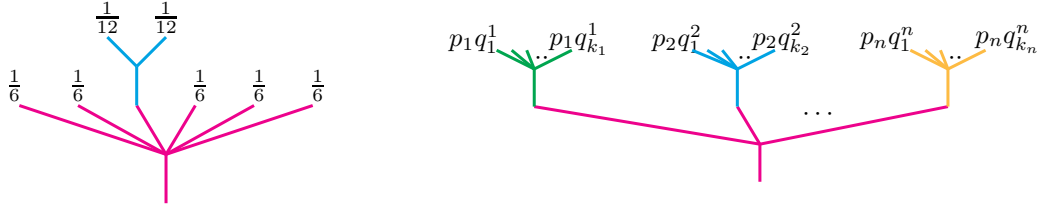
\noindent Further recall that if we have $n$ different distributions $q^i=(q^i_1,\ldots, q^i_{k_i})\in\Delta_{k_i}$ where $i=1,\ldots,n$, then we may compose them with $p$ simultaneously to obtain the following point in $\Delta_{k_1+\cdots+k_n},$

\[p\circ (q^1,\ldots,q^n)=(p_1q^1_1,\ldots,p_1q^1_{k_1},p_2q^2_1,\ldots,p_2q^2_{k_2},\ldots,p_nq^n_{1},\ldots,p_nq^n_{k_n}).\]
This simultaneous composition is illustrated by the tree on the right in Figure~\ref{fig:dice_pic}.
\end{Example}

Just as groups come to life when considering representations of them, so operads come to life when each abstract $n$-ary operation is mapped to a concrete $n$-ary operation on a particular object. This assignment is traditionally called an \textit{algebra} of the operad, but we prefer the more descriptive name \textit{representation.} 

\begin{Definition}Let $\O$ be an operad in {the category of sets}. A \emph{representation of $\O$}, or an \emph{$\O$-representation}, is {set} $X$ together with {functions}

\[\varphi_n\c \O(n)\to \End_X(n) \qquad \text{for $n\geq 1$}\]
that respect the operad unit and compositions. That is, $\varphi_n(1)=1$ and

\[\varphi_{n+m-1}(p\circ_i q)= \varphi_n(p)\circ_i\varphi_m(q)\]
for all $p\in\O(n), q\in\O(m)$ and $1\leq i \leq n$.
\end{Definition}
{Importantly, one may also wish to define a representation of an operad in any symmetric monoidal category $\mathsf{C}$ whenever ``$\End_X(n)$'' is in fact an object in $\mathsf{C}.$ It must consist of an object $X$ together with a family of morphisms $\mathcal{O}(n)\to\End_X(n)$ in $\mathsf{C}$ that are compatible with the operad unit and compositions. This holds, for instance, when the monoidal category $\mathsf{C}$ is also closed---that is, when it is equipped with an internal hom functor that is compatible with the monoidal product. Monoidal closure, however, will not be required in our work, which primarily concerns the category $\mathsf{Top}$ of topological spaces. Indeed, the} {main} example to have in mind is when $\O=\Delta$ is the operad of simplices and $X=\mathbb{R}$ is the real line in $\mathsf{Top}$. In this case, {we define} {$\End_\R(n):=\mathsf{Top}(\R^n,\R)$ to} be the space of continuous functions $\R^n\to \R$ equipped with the product topology. {Now, consider the} {continuous} maps $\varphi_n\c \Delta_n \to \End_\R(n)$ given by $p\mapsto \varphi_n(p)$ where $\varphi_n(p)(x):=\langle p,x\rangle= \sum_{i=1}^np_ix_i$ whenever $x=(x_1,\ldots,x_n)\in \R^n$. {Then, it is simple to check that $\varphi_{n+m-1}(p\circ_i q)= \varphi_n(p)\circ_i\varphi_m(q)$ for all $p,q,$ and $i$ and that $\varphi_n(1)=1$ for all $n$, and so $\mathbb{R}$ is a representation of $\Delta.$}

\section{Derivations of the Operad of Simplices}\label{sec:mainresult}

With these basic definitions in hand, the present goal is to define a mapping $d$ out of the topological operad $\Delta$ that satisfies an appropriate version of the Leibniz rule,
\begin{align}\label{eq:prodrule}
d(p\circ_iq)=dp\circ_i q + p\circ_idq \qquad \text{(desideratum)}
\end{align}
for all $p\in\Delta_n$ and $q\in\Delta_m$ {and for all $1\leq i \leq n$}. This desired equation suggests the codomain of $d$ should be a (bi)module over $\Delta$ that is, moreover, an abelian monoid. This motivates the following two definitions, the first of which is a slight generalization of that given by Markl in \cite{Mar}.

\begin{Definition} Let $\O=\{\mathcal{O}(1),\mathcal{O}(2),\ldots\}$ be an operad in a {symmetric monoidal} category $\C$. A \emph{{bi}module over $\O$}, or simply an \emph{$\O$-{bi}module}, is a collection of objects $M=\{M(1),M(2),\ldots\}$ in $\C$ together with morphisms 

 \begin{align*}
 \circ_i^L&= \O(n)\otimes M(m)\to M(n+m-1)\qquad \text{(left composition)}\\
 \circ_i^R&=M(n)\otimes \O(m)\to M(n+m-1)\qquad \text{(right composition)}
 \end{align*}
 in $\mathsf{C}$ for each $1\leq i\leq n$ {such} that whenever 
 \[
 p\otimes q\otimes r\in
 \begin{cases}
 M(n)\otimes \O(m)\otimes \O(k), \text{or}\\
 \O(n)\otimes M(m)\otimes\O(k),\text{or}\\
 \O(n)\otimes \O(m)\otimes M(k)
 \end{cases}
 \] 
the following holds:
 
 \begin{align}\label{eq:comp}
 (p\o_j q)\o_i r =
 \begin{cases}
 	(p\o_i r)\o_{j+k-1} q &\text{if $1\leq i \leq j-1$}\\
 	p\o_j(q \o_{i-j+1} r) &\text{if $j\leq i \leq j+m-1$}\\
 	(p\o_{i-m+1} r)\o_j q &\text{if $i \geq j+m$}.
 \end{cases}
 \end{align}
 \end{Definition}
The associativity requirements displayed in Equation (\ref{eq:comp})---and hence the intuition behind them---are completely analogous to those defining operads as illustrated in \mbox{Figure~\ref{fig:operad_axioms}}. The only difference here is that one of the three operations may come from the {bi}module rather than the operad. Here is the main example to have in mind.

\begin{Example}\label{ex:composition}
As every algebra is a bimodule over itself, so every representation of $\mathcal{O}$ is an $\mathcal{O}$-{bi}module in a straightforward way. Indeed, in the case of the topological operad of simplices, the maps comprising the $\Delta$-representation structure on $\mathbb{R}$ induce a $\Delta$-{bi}module structure on $\End_\R$. However, we will make use of a {slight variant of this} {bi}module structure. {Right composition will be defined in the expected way, though left composition will not. Explicitly,} we define the left and right composition maps

\begin{align*}
\o_i^L\c \Delta_n \times \mathsf{Top}(\R^m,\R)&\longrightarrow\mathsf{Top}(\R^{n+m-1},\R)\\
\o_i^R\c \mathsf{Top}(\R^n,\R) \times \Delta_m&\longrightarrow\mathsf{Top}(\R^{n+m-1},\R)
\end{align*}
as follows. Given a probability distribution $p\in\Delta_n$ and a continuous function $f\colon \mathbb{R}^m\to\mathbb{R}$, define left composition by {$p\circ_i^Lf:=\bar{p}\circ (0,\ldots,0,f,0,\ldots,0)$, where the composition on the right-hand side is defined as in the simultaneous composition in the endomorphism operad of $\mathbb{R}$ illustrated in Example} \ref{ex:endomorphism}{, and where each $0$ denotes the zero function $\mathbb{R}\to\mathbb{R}$. Here, recall that $\bar{p}\colon\mathbb{R}^n\to\mathbb{R}$ maps a point $x$ to the standard inner product $\langle p,x\rangle$ as introduced in Section} \ref{sec:intro}{. Unwinding this, left composition thus evaluates explicitly as} $(p\o_i^Lf)(x_1,\ldots,x_{n+m-1})=p_if(x_i,\ldots,x_{i+m-1})$. In words, the value of the left composite $p\circ_i^Lf\colon\mathbb{R}^{n+m-1}\to\mathbb{R}$ at a point $x$ is computed by evaluating $f$ at the $m$-subtuple of $x$ beginning at the $i^{\text{th}}$ coordinate and scaling that output by $p_i$. All other coordinates of $x$ are ignored. The picture to have in mind is that below, where the bold dots are imagined to be ``plugs'' that prevent the surplus coordinates from playing a role. In this picture, $n=3$ and $m=2$.

\begin{equation*}
(p\circ_2^L f)(x_1,x_2,x_3,x_4) \quad = \quad
\begin{tikzpicture}[x=.8cm,y=.6cm,baseline={(current bounding box.center)}]
	\node[] (a) at (0,0) {};
	\node[tinydot,color=magenta] (l) at (1,1) {};
	\node[tinydot,color=magenta] (r) at (-1,1) {};
	\node[] (m) at (0,1) {};
	\draw[thick, shorten <=-1.5mm, shorten >=-1mm, color=magenta] (0,-.75) -- (a);
	\draw[thick, shorten <=-1.5mm, color=magenta] (a) -- (l);
	\draw[thick, shorten <=-1.5mm, color=magenta] (a) -- (r);
	\draw[thick, shorten <=-1.5mm, color=magenta] (a) -- (m);

	\node[] (T) at (0,1.5) {};
	\draw[thick,shorten <=-1.5mm,shorten >=-1mm] (0,1) -- (T);
	\draw[thick,shorten <=-1.5mm] (T) -- (-.7,2);
	\draw[thick,shorten <=-1.5mm] (T) -- (.7,2);

	\node[left] at (0,-.8) {$p$}; 
	\node[left] at (0,1) {$f$}; 

	\node[above] at (-1,1) {{\color{gray}$x_1$}};
	\node[above] at (-.7,2) {$x_2$};
	\node[above] at (.7,2) {$x_3$};
	\node[above] at (1,1) {{\color{gray}$x_4$}};
\end{tikzpicture}
\quad = \quad p_2f(x_2,x_3) 
\end{equation*}

Given a probability distribution $q\in\Delta_m$ and a continuous function $g\colon \mathbb{R}^n\to\mathbb{R}$, define right composition by

\begin{align*}
(g\o_i^R q)(x_1,\ldots,&x_{n+m-1})\nonumber \\ 
&:=g(x_1,\ldots,x_{i-1},\sum_{k=1}^m q_kx_{i+k-1},x_{i+m},\ldots,x_{n+m-1}).
\end{align*}

{This} may be understood visually as well. The value of the right composite $g~\circ_i^R~q\colon\mathbb{R}^{n+m-1}\to\mathbb{R}$ at a point $x$ is computed by taking the inner product of $q$ with the $m$-tuple of $x$ beginning at the $i^{\text{th}}$ coordinate and using that number as the $i^{\text{th}}$ input of $g$ with all other coordinates of $x$ falling into place as in the picture below. There are no ``plugs'' in this instance since all coordinates of $x$ play a~role.

\begin{equation*}
(g\circ_2^R q)(x_1,x_2,x_3,x_4) \quad = \quad
\begin{tikzpicture}[x=.8cm,y=.6cm,baseline={(current bounding box.center)}]
	\node[] (a) at (0,0) {};
	\node[] (l) at (1,1) {};
	\node[] (r) at (-1,1) {};
	\node[] (m) at (0,1) {};
	\draw[thick, shorten <=-1.5mm, shorten >=-1mm] (0,-.75) -- (a);
	\draw[thick, shorten <=-1.5mm] (a) -- (l);
	\draw[thick, shorten <=-1.5mm] (a) -- (r);
	\draw[thick, shorten <=-1.5mm] (a) -- (m);

	\node[] (T) at (0,1.5) {};
	\draw[thick,shorten <=-1.5mm,shorten >=-1mm, color=cyan] (0,1) -- (T);
	\draw[thick,shorten <=-1.5mm, color=cyan] (T) -- (-.7,2);
	\draw[thick,shorten <=-1.5mm, color=cyan] (T) -- (.7,2);

	\node[left] at (0,-.8) {$g$}; 
	\node[left] at (0,1) {$q$}; 

	\node[above] at (-1,1) {$x_1$};
	\node[above] at (-.7,2) {$x_2$};
	\node[above] at (.7,2) {$x_3$};
	\node[above] at (1,1) {$x_4$};
\end{tikzpicture}
\quad=\quad g(x_1,q_1x_2+q_2x_3,x_4)
\end{equation*}

These examples suggest the inner product notation is a convenient choice. Given $N\geq 1$ and $k\leq N$ and a point $x\in\mathbb{R}^N$, let $\mathbf{x}_{i,k}\in\mathbb{R}^k$ denote the $k$-subtuple of $x$ beginning at the $i^{\text{th}}$ coordinate:

\[\mathbf{x}_{i,k}:=(x_i,\ldots,x_{i+k-1}).\]
Then given any point $x\in\mathbb{R}^{n+m-1}$, the left and right composition maps may be written more succinctly as 

\begin{align*}
(p\circ_i^Lf)(x) &= p_if(\mathbf{x}_{i,m})\\
(g\circ_i^R q)(x) &= g(x_1,\ldots,x_{i-1},\langle q,\mathbf{x}_{i,m}\rangle,x_{i+m},\ldots,x_{n+m-1}).
\end{align*}
We will use this notation below and will always write $\mathbf{x}_{i}$ in lieu of $\mathbf{x}_{i,m}$ since the context will make it clear that $\mathbf{x}_{i}$ must be an $m$-tuple. The boldface font is used to distinguish a tuple $\mathbf{x}_i$ from a real number $x_i$. Finally, note that the maps $\circ_i^L$ and $\circ_i^R$ are continuous since $f$ and $g$ are continuous, and moreover that the {associativity} requirements in Equation (\ref{eq:comp}) are analogous to those illustrated in Figure~\ref{fig:operad_axioms}, so it is straightforward to verify they are satisfied. {In particular, the zero functions appearing in the definition of $\circ_i^L$ simplify the situation greatly. For instance, several of the associativity requirements follow from the simple fact that multiplying an input $x_i$ by a probability and then mapping the result to zero is the same as first mapping the input to zero and then multiplying that zero by a probability.} So $\End_\R$ is indeed a $\Delta$-{bi}module.
\end{Example}

Next, recall that the desired Leibniz rule in Equation (\ref{eq:prodrule}) suggests the {bi}module should be equipped with a notion of addition. This motivates the following definition.

\begin{Definition} Let $\O$ be an operad in a {symmetric monoidal} category $\C$. An $\O$-{bi}module $M$ is an \emph{abelian $\O$-{bi}module}
if each $M(n)$ is an abelian monoid in $\C$; that is, if for each $n=1,2,\ldots$ the following hold:
	\begin{itemize}
	\item[(i)] [associativity, commutativity] there is a morphism $\mu_n\c M(n)\times M(n)\to M(n)$ in $\C$ such that $\mu_n(\mu_n(a,b),c)=\mu_n(a,\mu_n(b,c))$ and $\mu_n(a,b)=\mu_n(b,a)$ for all $a,b,c\in M(n)$,
	\item[(ii)] [identity] there is an element $1\in M(n)$ {such} that $\mu_n(1,a)=a=\mu_n(a,1)$ for all $a\in M(n)$.
	\end{itemize}
\end{Definition}
As the primary example, consider $\End_\R$ viewed as a $\Delta$-{bi}module as described in Example \ref{ex:composition}. For each $n,$ define $\mu_n\c \End_\R(n)\times \End_\R(n)\to \End_\R(n)$ by pointwise addition, meaning that for each  $f,g\in \End_\R(n)$ we have $\mu_n(f,g)=f+g$ where $(f+g)(x):=f(x)+ g(x)$ for all $x\in \R^n$. The identity element in $\End_\R(n)$ is the constant map at zero. Moreover each $\mu_n$ is continuous and inherits associativity and commutativity from $\R.$ In this way, $\End_\R$ is an abelian $\Delta$-{bi}module. 

\begin{Remark}\label{rem:distribute}
Notice that the $\Delta$-{bi}module composition maps $\circ_i^L$ and $\circ_i^R$ distribute over sums in the abelian $\Delta$-{bi}module $\End_\R$. In other words, for all continuous functions $f,g\in\End_\R(n)$ and for all probability distributions $q\in\Delta_m$,

\begin{align*}
(f+g)\circ_i^R q=f\circ_i^Rq + g\circ_i^R q, \qquad1\leq i \leq n
\end{align*}
and similarly for left composition $\circ_i^L$. This follows directly from pointwise addition. 
\end{Remark}

With this setup in mind, our desideratum in Equation (\ref{eq:prodrule}) is now realized in the following definition.

\begin{Definition}
Let $\O$ be an operad in a category $\C$ and let $M$ be an abelian $\O$-{bi}module. A \emph{derivation of $\O$ valued in $M$} is sequence of morphisms  $\{d_n\c \O(n)\to M(n)\}$ in $\C$ satisfying

\begin{align}\label{eq:deriv}
d_{n+m-1}(p\circ_i q)= d_np\circ_i^Rq + p\circ_i^L d_mq
\end{align}
for all $p\in \O(n), q\in \O(m)$ and for all $1\leq i \leq n.$
\end{Definition}
In the special case when $\O$ is a linear operad, this definition coincides with that given by Markl in \cite{Mar}.
In what follows, we omit the subscripts and simply write $d$ instead of $d_n$. Now, suppose $\mathcal{O}=\Delta$ is the operad of topological simplices and {$\End_\mathbb{R}$} is equipped with the structure of an abelian $\Delta$-{bi}module given above. Here is the picture to have in mind for Equation (\ref{eq:deriv}):

\begin{equation*}
d\left(
\begin{tikzpicture}[x=.6cm,y=.35cm,baseline={(current bounding box.center)}]
	\node[] (a) at (0,0) {};
	\node[] (l) at (1,1) {};
	\node[] (r) at (-1,1) {};
	\node[] (m) at (0,1) {};
	\draw[thick, shorten <=-1.5mm, shorten >=-1mm, color=magenta] (0,-.75) -- (a);
	\draw[thick, shorten <=-1.5mm, color=magenta] (a) -- (l);
	\draw[thick, shorten <=-1.5mm, color=magenta] (a) -- (r);
	\draw[thick, shorten <=-1.5mm, color=magenta] (a) -- (m);

	\node[] (T) at (0,1.5) {};
	\draw[thick,color=Cerulean,shorten <=-1.5mm,shorten >=-1mm] (0,1) -- (T);
	\draw[thick,color=Cerulean,shorten <=-1.5mm] (T) -- (-.7,2);
	\draw[thick,color=Cerulean,shorten <=-1.5mm] (T) -- (.7,2);
\end{tikzpicture}
\right)
\quad = \quad
\begin{tikzpicture}[x=.6cm,y=.35cm,baseline={(current bounding box.center)}]
	\node[] (a) at (0,0) {};
	\node[] (l) at (1,1) {};
	\node[] (r) at (-1,1) {};
	\node[] (m) at (0,1) {};
	\draw[thick, shorten <=-1.5mm, shorten >=-1mm, color=magenta] (0,-.75) -- (a);
	\draw[thick, shorten <=-1.5mm, color=magenta] (a) -- (l);
	\draw[thick, shorten <=-1.5mm, color=magenta] (a) -- (r);
	\draw[thick, shorten <=-1.5mm, color=magenta] (a) -- (m);

	\node[] (T) at (0,1.5) {};
	\draw[thick,color=Cerulean,shorten <=-1.5mm,shorten >=-1mm] (0,1) -- (T);
	\draw[thick,color=Cerulean,shorten <=-1.5mm] (T) -- (-.7,2);
	\draw[thick,color=Cerulean,shorten <=-1.5mm] (T) -- (.7,2);

	\node [] at (-.5,-.75) {$d$};
\end{tikzpicture}
\quad + \quad
\begin{tikzpicture}[x=.6cm,y=.35cm,baseline={(current bounding box.center)}]
	\node[] (a) at (0,0) {};
	\node[tinydot,fill=magenta,inner sep=.6mm,] (l) at (1,1) {};
	\node[tinydot,fill=magenta,inner sep=.6mm,] (r) at (-1,1) {};
	\node[] (m) at (0,1) {};
	\draw[thick, shorten <=-1.5mm, shorten >=-1mm, color=magenta] (0,-.75) -- (a);
	\draw[thick, shorten <=-1.5mm, color=magenta] (a) -- (l);
	\draw[thick, shorten <=-1.5mm, color=magenta] (a) -- (r);
	\draw[thick, shorten <=-1.5mm, color=magenta] (a) -- (m);

	\node[] (T) at (0,1.5) {};
	\draw[thick,color=Cerulean,shorten <=-1.5mm,shorten >=-1mm] (0,1) -- (T);
	\draw[thick,color=Cerulean,shorten <=-1.5mm] (T) -- (-.7,2);
	\draw[thick,color=Cerulean,shorten <=-1.5mm] (T) -- (.7,2);

	\node [] at (-.5,1.2) {$d$};
\end{tikzpicture}
\end{equation*}
On the right-hand side we have used the ``plug'' notation introduced in Example \ref{ex:composition}, which can also be understood explicitly by evaluating $d$ at a point $x\in\mathbb{R}^{n+m-1}$,

\begin{align*}
d(p\circ_i q)(x) &= (dp\circ_i^Rq)(x) + (p\circ_i^L dq)(x)\\
&=dp(x_1,\ldots,\langle q,\mathbf{x}_i\rangle,\ldots,x_{n+m-1}) + p_idq(\mathbf{x}_i).
\end{align*}

Of particular interest is the behavior of a derivation $\{d\colon\Delta_n\to\End_{\mathbb{R}}(n)\}$ when it is applied to a simultaneous composition of probability distributions. A derivation applied to the composite $(p\circ_j q)\circ_i r$ for probability distributions $p\in\Delta_n,q\in\Delta_m$, and $r\in\Delta_k$ can be understood in a convenient picture when $q$ and $r$ are composed onto different leaves of $p$; that is, when $1\leq i\leq j-1$ or $i\geq j+m$. This follows straightforwardly from a repeated application of $d$. Indeed, by definition we have $d((p\circ_j q)\circ_i r)=d(p\circ_j q)\circ_i^R r + (p\circ_j q) \circ_i^L dr$ and by applying the Leibniz rule again to the first summand, this is equal to $(dp\circ_j^R q + p\circ_j^Ldq)\circ_i^R r + (p\circ_j q) \circ_i^L dr,$ which we can expand to obtain $(dp\circ_j^R q)\circ_i^R r + (p\circ_j^Ldq)\circ_i^R r + (p\circ_j q) \circ_i^L dr$ since composition distributes over sums as noted in \mbox{Remark \ref{rem:distribute}}. We will identify this function with the picture below in lieu of the cumbersome notation.

\begin{equation*}
d\left(
\begin{tikzpicture}[x=.55cm,y=.6cm,baseline={(current bounding box.center)},scale=0.6]
	\node[] (a) at (0,0) {};
	\node[] (r) at (2.5,1) {};
	\node[] (l) at (-2.5,1) {};
	\node[] (m) at (0,1) {};
	\draw[thick, shorten <=-1.5mm, shorten >=-1mm, color=magenta] (0,-.75) -- (a);
	\draw[thick, shorten <=-1.5mm, color=magenta] (a) -- (r);
	\draw[thick, shorten <=-1.5mm, color=magenta] (a) -- (l);
	\draw[thick, shorten <=-1.5mm, color=magenta] (a) -- (m);
	
	\node[] (T) at (0,1.5) {};
	\draw[thick,color=Cerulean,shorten <=-1.5mm,shorten >=-1.2mm] (0,1) -- (T);
	\draw[thick,color=Cerulean,shorten <=-1.5mm] (T) -- (-.7,2);
	\draw[thick,color=Cerulean,shorten <=-1.5mm] (T) -- (.7,2);

	\node[] (p1) at (-2.1,1.5) {};
	\draw[thick,color=Green,shorten <=-1.1mm,shorten >=-1.3mm] (-2.1,1.2) -- (p1);
	\draw[thick,color=Green,shorten <=-1.5mm] (p1) -- (-2.8,2);
	\draw[thick,color=Green,shorten <=-1.5mm] (p1) -- (-1.4,2);
\end{tikzpicture}
\right)
\quad = \quad
\begin{tikzpicture}[x=.55cm,y=.6cm,baseline={(current bounding box.center)},scale=0.6]
	\node[] (a) at (0,0) {};
	\node[] (r) at (2.5,1) {};
	\node[] (l) at (-2.5,1) {};
	\node[] (m) at (0,1) {};
	\draw[thick, shorten <=-1.5mm, shorten >=-1mm, color=magenta] (0,-.75) -- (a);
	\draw[thick, shorten <=-1.5mm, color=magenta] (a) -- (r);
	\draw[thick, shorten <=-1.5mm, color=magenta] (a) -- (l);
	\draw[thick, shorten <=-1.5mm, color=magenta] (a) -- (m);

	\node[] (T) at (0,1.5) {};
	\draw[thick,color=Cerulean,shorten <=-1.5mm,shorten >=-1.2mm] (0,1) -- (T);
	\draw[thick,color=Cerulean,shorten <=-1.5mm] (T) -- (-.7,2);
	\draw[thick,color=Cerulean,shorten <=-1.5mm] (T) -- (.7,2);

	\node[] (p1) at (-2.1,1.5) {};
	\draw[thick,color=Green,shorten <=-1.1mm,shorten >=-1.3mm] (-2.1,1.2) -- (p1);
	\draw[thick,color=Green,shorten <=-1.5mm] (p1) -- (-2.8,2);
	\draw[thick,color=Green,shorten <=-1.5mm] (p1) -- (-1.4,2);

	\node [] at (-.8,-.8) {$d$};
\end{tikzpicture}
\quad + \quad
\begin{tikzpicture}[x=.55cm,y=.6cm,baseline={(current bounding box.center)},scale=0.6]
	\node[] (a) at (0,0) {};
	\node[tinydot,inner sep=.6mm,fill=magenta] (r) at (2.5,1) {};
	\node[] (l) at (-2.5,1) {};
	\node[] (m) at (0,1) {};
	\draw[thick, shorten <=-1.5mm, shorten >=-1mm, color=magenta] (0,-.75) -- (a);
	\draw[thick, shorten <=-1.5mm, color=magenta] (a) -- (r);
	\draw[thick, shorten <=-1.5mm, color=magenta] (a) -- (l);
	\draw[thick, shorten <=-1.5mm, color=magenta] (a) -- (m);

	\node[tinydot,inner sep=.6mm,fill=Cerulean] at (-.7,2) {};
	\node[tinydot,inner sep=.6mm,fill=Cerulean] at (.7,2) {};

	\node[] (T) at (0,1.5) {};
	\draw[thick,color=Cerulean,shorten <=-1.5mm,shorten >=-1.2mm] (0,1) -- (T);
	\draw[thick,color=Cerulean,shorten <=-1.5mm] (T) -- (-.7,2);
	\draw[thick,color=Cerulean,shorten <=-1.5mm] (T) -- (.7,2);

	\node[] (p1) at (-2.1,1.5) {};
	\draw[thick,color=Green,shorten <=-1.1mm,shorten >=-1.3mm] (-2.1,1.2) -- (p1);
	\draw[thick,color=Green,shorten <=-1.5mm] (p1) -- (-2.8,2);
	\draw[thick,color=Green,shorten <=-1.5mm] (p1) -- (-1.4,2);

	\node [] at (-3,1.3) {$d$};
\end{tikzpicture}
\quad + \quad
\begin{tikzpicture}[x=.55cm,y=.6cm,baseline={(current bounding box.center)},scale=0.6]
	\node[] (a) at (0,0) {};
	\node[tinydot,inner sep=.6mm,fill=magenta] (r) at (2.5,1) {};
	\node[] (l) at (-2.5,1) {};
	\node[] (m) at (0,1) {};
	\draw[thick, shorten <=-1.5mm, shorten >=-1mm, color=magenta] (0,-.75) -- (a);
	\draw[thick, shorten <=-1.5mm, color=magenta] (a) -- (r);
	\draw[thick, shorten <=-1.5mm, color=magenta] (a) -- (l);
	\draw[thick, shorten <=-1.5mm, color=magenta] (a) -- (m);

	\node[] (T) at (0,1.5) {};
	\draw[thick,color=Cerulean,shorten <=-1.5mm,shorten >=-1.2mm] (0,1) -- (T);
	\draw[thick,color=Cerulean,shorten <=-1.5mm] (T) -- (-.7,2);
	\draw[thick,color=Cerulean,shorten <=-1.5mm] (T) -- (.7,2);

	\node[] (p1) at (-2.1,1.5) {};
	\draw[thick,color=Green,shorten <=-1.1mm,shorten >=-1.3mm] (-2.1,1.2) -- (p1);
	\draw[thick,color=Green,shorten <=-1.5mm] (p1) -- (-2.8,2);
	\draw[thick,color=Green,shorten <=-1.5mm] (p1) -- (-1.4,2);

	\node[tinydot,inner sep=.6mm,fill=Green] at (-2.8,2) {};
	\node[tinydot,inner sep=.6mm,fill=Green] at (-1.4,2) {};

	\node [] at (-.8,1.3) {$d$};
\end{tikzpicture}
\end{equation*}
Importantly, the {obvious generalization of the} formula holds for any simultaneous composition $p\circ(q^1,\ldots,q^n)$ for any $p\in\Delta_n$ and $q^i\in\Delta_{k_i}$ where $i=1,\ldots,n$. This again follows directly from repeated applications of Equation (\ref{eq:deriv}), as illustrated below.

\begin{equation*}
d\left(
\begin{tikzpicture}[x=.55cm,y=.6cm,baseline={(current bounding box.center)},scale=0.5]
	\node[] (a) at (0,0) {};
	\node[] (r) at (2.5,1) {};
	\node[] (l) at (-2.5,1) {};
	\node[] (m) at (0,1) {};
	\draw[thick, shorten <=-1.5mm, shorten >=-1mm, color=magenta] (0,-.75) -- (a);
	\draw[thick, shorten <=-1.5mm, color=magenta] (a) -- (r);
	\draw[thick, shorten <=-1.5mm, color=magenta] (a) -- (l);
	\draw[thick, shorten <=-1.5mm, color=magenta] (a) -- (m);
	
	\node[] (T) at (0,1.5) {};
	\draw[thick,color=Cerulean,shorten <=-1.5mm,shorten >=-1.2mm] (0,1) -- (T);
	\draw[thick,color=Cerulean,shorten <=-1.5mm] (T) -- (-.7,2);
	\draw[thick,color=Cerulean,shorten <=-1.5mm] (T) -- (.7,2);

	\node[] (p1) at (-2,1.5) {};
	\draw[thick,color=Green,shorten <=-1.1mm,shorten >=-1.3mm] (-2,1.2) -- (p1);
	\draw[thick,color=Green,shorten <=-1.5mm] (p1) -- (-2.8,2);
	\draw[thick,color=Green,shorten <=-1.5mm] (p1) -- (-1.4,2);

	\node[] (p3) at (2,1.5) {};
	\draw[thick,color=Dandelion,shorten <=-1.1mm,shorten >=-1.3mm] (2,1.2) -- (p3);
	\draw[thick,color=Dandelion,shorten <=-1.5mm] (p3) -- (2.8,2);
	\draw[thick,color=Dandelion,shorten <=-1.5mm] (p3) -- (1.4,2);
\end{tikzpicture}
\right)
\quad = \quad
\begin{tikzpicture}[x=.55cm,y=.6cm,baseline={(current bounding box.center)},scale=0.5]
	\node[] (a) at (0,0) {};
	\node[] (r) at (2.5,1) {};
	\node[] (l) at (-2.5,1) {};
	\node[] (m) at (0,1) {};
	\draw[thick, shorten <=-1.5mm, shorten >=-1mm, color=magenta] (0,-.75) -- (a);
	\draw[thick, shorten <=-1.5mm, color=magenta] (a) -- (r);
	\draw[thick, shorten <=-1.5mm, color=magenta] (a) -- (l);
	\draw[thick, shorten <=-1.5mm, color=magenta] (a) -- (m);

	\node[] (T) at (0,1.5) {};
	\draw[thick,color=Cerulean,shorten <=-1.5mm,shorten >=-1.2mm] (0,1) -- (T);
	\draw[thick,color=Cerulean,shorten <=-1.5mm] (T) -- (-.7,2);
	\draw[thick,color=Cerulean,shorten <=-1.5mm] (T) -- (.7,2);

	\node[] (p1) at (-2,1.5) {};
	\draw[thick,color=Green,shorten <=-1.1mm,shorten >=-1.3mm] (-2,1.2) -- (p1);
	\draw[thick,color=Green,shorten <=-1.5mm] (p1) -- (-2.8,2);
	\draw[thick,color=Green,shorten <=-1.5mm] (p1) -- (-1.4,2);

	\node[] (p3) at (2,1.5) {};
	\draw[thick,color=Dandelion,shorten <=-1.1mm,shorten >=-1.3mm] (2,1.2) -- (p3);
	\draw[thick,color=Dandelion,shorten <=-1.5mm] (p3) -- (2.8,2);
	\draw[thick,color=Dandelion,shorten <=-1.5mm] (p3) -- (1.4,2);

	\node [] at (-.8,-.8) {$d$};
\end{tikzpicture}
\quad + \quad
\begin{tikzpicture}[x=.55cm,y=.6cm,baseline={(current bounding box.center)},scale=0.5]
	\node[] (a) at (0,0) {};
	\node[] (r) at (2.5,1) {};
	\node[] (l) at (-2.5,1) {};
	\node[] (m) at (0,1) {};
	\draw[thick, shorten <=-1.5mm, shorten >=-1mm, color=magenta] (0,-.75) -- (a);
	\draw[thick, shorten <=-1.5mm, color=magenta] (a) -- (r);
	\draw[thick, shorten <=-1.5mm, color=magenta] (a) -- (l);
	\draw[thick, shorten <=-1.5mm, color=magenta] (a) -- (m);

	\node[] (T) at (0,1.5) {};
	\draw[thick,color=Cerulean,shorten <=-1.5mm,shorten >=-1.2mm] (0,1) -- (T);
	\draw[thick,color=Cerulean,shorten <=-1.5mm] (T) -- (-.7,2);
	\draw[thick,color=Cerulean,shorten <=-1.5mm] (T) -- (.7,2);
	\node[tinydot,fill=Cerulean,inner sep = .5mm] at (-.7,2) {};
	\node[tinydot,fill=Cerulean,inner sep = .5mm] at (.7,2) {};

	\node[] (p1) at (-2,1.5) {};
	\draw[thick,color=Green,shorten <=-1.1mm,shorten >=-1.3mm] (-2,1.2) -- (p1);
	\draw[thick,color=Green,shorten <=-1.5mm] (p1) -- (-2.8,2);
	\draw[thick,color=Green,shorten <=-1.5mm] (p1) -- (-1.4,2);

	\node[] (p3) at (2,1.5) {};
	\draw[thick,color=Dandelion,shorten <=-1.1mm,shorten >=-1.3mm] (2,1.2) -- (p3);
	\draw[thick,color=Dandelion,shorten <=-1.5mm] (p3) -- (2.8,2);
	\draw[thick,color=Dandelion,shorten <=-1.5mm] (p3) -- (1.4,2);
	\node[tinydot,fill=Dandelion,inner sep = .5mm] at (2.8,2) {};
	\node[tinydot,fill=Dandelion,inner sep = .5mm] at (1.4,2) {};

	\node [] at (-3,1.3) {$d$};
\end{tikzpicture}
\quad + \quad
\begin{tikzpicture}[x=.55cm,y=.6cm,baseline={(current bounding box.center)},scale=0.5]
	\node[] (a) at (0,0) {};
	\node[] (r) at (2.5,1) {};
	\node[] (l) at (-2.5,1) {};
	\node[] (m) at (0,1) {};
	\draw[thick, shorten <=-1.5mm, shorten >=-1mm, color=magenta] (0,-.75) -- (a);
	\draw[thick, shorten <=-1.5mm, color=magenta] (a) -- (r);
	\draw[thick, shorten <=-1.5mm, color=magenta] (a) -- (l);
	\draw[thick, shorten <=-1.5mm, color=magenta] (a) -- (m);

	\node[] (T) at (0,1.5) {};
	\draw[thick,color=Cerulean,shorten <=-1.5mm,shorten >=-1.2mm] (0,1) -- (T);
	\draw[thick,color=Cerulean,shorten <=-1.5mm] (T) -- (-.7,2);
	\draw[thick,color=Cerulean,shorten <=-1.5mm] (T) -- (.7,2);

	\node[] (p1) at (-2,1.5) {};
	\draw[thick,color=Green,shorten <=-1.1mm,shorten >=-1.3mm] (-2,1.2) -- (p1);
	\draw[thick,color=Green,shorten <=-1.5mm] (p1) -- (-2.8,2);
	\draw[thick,color=Green,shorten <=-1.5mm] (p1) -- (-1.4,2);
	\node[tinydot,fill=Green,inner sep = .5mm] at (-2.8,2) {};
	\node[tinydot,fill=Green,inner sep = .5mm] at (-1.4,2) {};

	\node[] (p3) at (2,1.5) {};
	\draw[thick,color=Dandelion,shorten <=-1.1mm,shorten >=-1.3mm] (2,1.2) -- (p3);
	\draw[thick,color=Dandelion,shorten <=-1.5mm] (p3) -- (2.8,2);
	\draw[thick,color=Dandelion,shorten <=-1.5mm] (p3) -- (1.4,2);
	\node[tinydot,fill=Dandelion,inner sep = .5mm] at (2.8,2) {};
	\node[tinydot,fill=Dandelion,inner sep = .5mm] at (1.4,2) {};

	\node [] at (-.8,1.3) {$d$};
\end{tikzpicture}
\quad + \quad
\begin{tikzpicture}[x=.55cm,y=.6cm,baseline={(current bounding box.center)},scale=0.5]
	\node[] (a) at (0,0) {};
	\node[] (r) at (2.5,1) {};
	\node[] (l) at (-2.5,1) {};
	\node[] (m) at (0,1) {};
	\draw[thick, shorten <=-1.5mm, shorten >=-1mm, color=magenta] (0,-.75) -- (a);
	\draw[thick, shorten <=-1.5mm, color=magenta] (a) -- (r);
	\draw[thick, shorten <=-1.5mm, color=magenta] (a) -- (l);
	\draw[thick, shorten <=-1.5mm, color=magenta] (a) -- (m);

	\node[] (T) at (0,1.5) {};
	\draw[thick,color=Cerulean,shorten <=-1.5mm,shorten >=-1.2mm] (0,1) -- (T);
	\draw[thick,color=Cerulean,shorten <=-1.5mm] (T) -- (-.7,2);
	\draw[thick,color=Cerulean,shorten <=-1.5mm] (T) -- (.7,2);
	\node[tinydot,fill=Cerulean,inner sep = .5mm] at (-.7,2) {};
	\node[tinydot,fill=Cerulean,inner sep = .5mm] at (.7,2) {};

	\node[] (p1) at (-2,1.5) {};
	\draw[thick,color=Green,shorten <=-1.1mm,shorten >=-1.3 mm] (-2,1.2) -- (p1);
	\draw[thick,color=Green,shorten <=-1.5mm] (p1) -- (-2.8,2);
	\draw[thick,color=Green,shorten <=-1.5mm] (p1) -- (-1.4,2);
	\node[tinydot,fill=Green,inner sep = .5mm] at (-2.8,2) {};
	\node[tinydot,fill=Green,inner sep = .5mm] at (-1.4,2) {};

	\node[] (p3) at (2,1.5) {};
	\draw[thick,color=Dandelion,shorten <=-1.1mm,shorten >=-1.3mm] (2,1.2) -- (p3);
	\draw[thick,color=Dandelion,shorten <=-1.5mm] (p3) -- (2.8,2);
	\draw[thick,color=Dandelion,shorten <=-1.5mm] (p3) -- (1.4,2);

	\node [] at (1,1.3) {$d$};
\end{tikzpicture}
\end{equation*}
This is summarized in the following proposition.

\begin{Proposition}\label{proposition}
Let $p\in \Delta_n$ and $q^i\in\Delta_{k_i}$ for $n,k_1,\ldots,k_n\geq 1$ and let $\{d\colon \Delta_n\to\End_\mathbb{R}(n)\}$ be a derivation of the operad of topological simplices. Then for any point $x\in\mathbb{R}^{k_1+\cdots+k_n}$,

\begin{equation*}
d(p\circ (q^1,\ldots,q^n))(x)=dp(\langle q^1,\mathbf{x}_1\rangle,\cdots,\langle q^n,\mathbf{x}_n\rangle) + \sum_{i=1}^np_idq^i(\mathbf{x}_i).
\end{equation*}

\end{Proposition}
\noindent Finally, the main result follows.

\begin{Theorem}\label{thm:main2}
Shannon entropy defines a derivation of the operad of topological simplices, and for every derivation of this operad there exists a point at which it is given by a constant multiple of Shannon entropy.
\end{Theorem}

\begin{proof}
For each $n\geq 1$ define $d\colon \Delta_n\to\End_\mathbb{R}(n)$ by $p\mapsto dp$ where $dp(x)=H(p)$ is constant for all $x\in\mathbb{R}^n$. Then, $d$ is continuous since $H$ is continuous. 
Moreover, if $p=(p_1,\ldots,p_n)\in \Delta_n$ and $q=(q_1,\ldots,q_m)\in \Delta_m$ are probability distributions, then for any $x\in \R^{m+n-1}$ and $1\leq i \leq n$, we have

\begin{align*}
d(p\o_i q)(x)= H(p\o_i q) &=-\left( \sum_{k=1}^{i-1}p_k\log p_k+p_i\sum_{k=1}^mq_k\log(p_iq_k) + \sum_{k=i+1}^{n}p_k\log p_k\right)\\[10pt]
&=-\left( \sum_{k=1}^{i-1}p_k\log p_k +  p_i\log p_i\sum_{k=1}^mq_k + p_i\sum_{k=1}^mq_k\log q_k  + \sum_{k=i+1}^{n}p_k\log p_k\right)\\[10pt]
&=-\left(\sum_{k=1}^n p_k\log p_k + p_i\sum_{k=1}^mq_k\log q_k  \right)\\[10pt]
&= H(p) + p_iH(q)\\[10pt]
&= (dp\o_i^R q + p\o_i^Ldq)(x),
\end{align*}
where the last line follows since $(dp\o_i^Rq)(x)$ is computed by evaluating the function $dp$ at some point, and this function is assumed to be constant at $H(p).$

Conversely, suppose $\{d\c \Delta_n\to \End_\R(n)\}$ is a derivation. For each $n\geq 1$ define a function $F\colon\Delta_n\to\mathbb{R}$ by $F(p)=dp(0)$ where $0=(0,\ldots,0)\in\mathbb{R}^n$. Then $F$ is continuous since $d$ is continuous, and Proposition \ref{proposition} further implies that

\begin{align*}
F(p\circ (q^1,\ldots,q^n)) &= d(p\circ (q^1,\ldots,q^n))(0)\\
&= dp(\langle q^1,\mathbf{0}_1\rangle,\ldots,\langle q^n,\mathbf{0}_n\rangle) + \sum_{i=1}^np_idq^i(\mathbf{0}_i)\\
&=dp(0) + \sum_{i=1}^np_idq^i(0)\\
&=F(p) + \sum_{i=1}^np_iF(q^i).
\end{align*}

From the Faddeev--Leinster result in Theorem \ref{thm:FL}, it follows that $dp(0)=F(p)=cH(p)$ for some $c\in\mathbb{R}$. 
\end{proof}

Notice that the important Equation (\ref{eq:Baez}) mentioned in the introduction is obtained as a corollary. Indeed, if for each $n\geq 1$ the map $d\colon\Delta_n\to\End_\mathbb{R}(n)$ is defined to be constant at entropy $p\mapsto dp\equiv H(p)$, then $d$ is a derivation by Theorem \ref{thm:main2} and so Proposition \ref{proposition} yields the following by evaluating $d(p\circ (q^1,\ldots,q^n))$ at any point.
\begin{Corollary*}Let $p\in\Delta_n$ and $q^i\in \Delta_{k_i}$ with $1\leq i \leq n.$ Then

\[H(p\circ (q^1,\ldots,q^n))=H(p)+ \sum_{i=1}^np_iH(q^i).\]

\end{Corollary*}

As a closing remark, Faddeev's characterization of entropy in Theorem~\ref{thm:FL} can be reexpressed using the language of category theory and operads as in \cite{leinster2021entropy}, (Theorem 12.3.1). We have omitted this language here but invite the reader to explore the full category theoretical story in Chapter 12 of Leinster's book.

\bibliographystyle{alpha}
\bibliography{references}{}

\begin{thebibliography}{EVG15}

\bibitem[Bae11]{EFunctor}
John~C. Baez.
\newblock Entropy as a functor, 2011.
\newblock Blog post. Available online:
  \url{https://www.ncatlab.org/johnbaez/show/Entropy+as+a+functor}.

\bibitem[BB15]{baudot2015}
Pierre Baudot and Daniel Bennequin.
\newblock The homological nature of entropy.
\newblock {\em Entropy}, 17(5):3253--3318, 2015.

\bibitem[BFL11]{BFL}
John~C. Baez, Tobias Fritz, and Tom Leinster.
\newblock A characterization of entropy in terms of information loss.
\newblock {\em Entropy}, 13:1945--1957, 2011.
\newblock doi:10.3390/e13111945.

\bibitem[BV73]{BV73}
J.~M. Boardman and R.~Vogt.
\newblock Homotopy invariant algebraic structures on topological spaces.
\newblock volume 347 of {\em Lecture Notes in Mathematics}. Springer, 1973.

\bibitem[EVG15]{poly2015}
Philippe Elbaz-Vincent and Herbet Gangl.
\newblock Finite polylogarithms, their multiple analogues and the {S}hannon
  entropy.
\newblock In F.~Nielsen and F.~Barbaresco, editors, {\em Geometric Science of
  Information. GSI 2015.}, volume 9389 of {\em Lecture Notes in Computer
  Science}, pages 277--285. Springer, Cham., 2015.

\bibitem[Fad56]{Fad}
D.~K. Faddeev.
\newblock On the concept of entropy of a finite probabilistic scheme.
\newblock {\em Uspekhi Mat. Nauk}, 11:227--231, 1956.
\newblock (In Russian).

\bibitem[Lam69]{lambek}
Joachim Lambek.
\newblock Deductive systems and categories ii. standard constructions and
  closed categories.
\newblock In P.~Hilton, editor, {\em Category Theory, Homology Theory and their
  Applications, I (Battelle Institute Conference, Seattle, 1968)}, volume~68 of
  {\em Lecture Notes in Mathematics}. Springer, 1969.

\bibitem[Lei21]{leinster2021entropy}
Tom Leinster.
\newblock {\em Entropy and Diversity: The Axiomatic Approach}.
\newblock Cambridge University Press, 2021.

\bibitem[LV12]{loday2012algebraic}
Jean-Louis Loday and Bruno Vallette.
\newblock {\em Algebraic Operads}.
\newblock Grundlehren der mathematischen Wissenschaften. Springer Berlin
  Heidelberg, 2012.

\bibitem[Mai19]{Mainiero}
Tom Mainiero.
\newblock Homological tools for the quantum mechanic.
\newblock 2019.
\newblock arXiv:1901.02011.

\bibitem[Mar96]{Mar}
Martin Markl.
\newblock Models for operads.
\newblock {\em Communications in Algebra}, 24(4):1471--1500, 1996.
\newblock arXiv: \url{arxiv.org/abs/hep-th/9411208}.

\bibitem[May72]{may2006geometry}
J.P. May.
\newblock The geometry of iterated loop spaces.
\newblock volume 271 of {\em Lecture Notes in Mathematics}. Springer, 1972.

\bibitem[MSS02]{markl2002operads}
Martin Markl, Steven Shnider, and Jim Stasheff.
\newblock {\em Operads in Algebra, Topology and Physics}.
\newblock Mathematical surveys and monographs. American Mathematical Society,
  2002.

\bibitem[Par20]{arthur}
Arthur~J. Parzygnat.
\newblock A functorial characterization of von neumann entropy.
\newblock 2020.
\newblock arXiv:2009.07125.

\bibitem[Sta04]{Stasheff}
Jim Stasheff.
\newblock What is... an operad?
\newblock {\em Notices Amer. Math. Soc.}, 51:630--631, 2004.

\bibitem[Val12]{Bruno}
Bruno Vallette.
\newblock Algebra + homotopy = operad.
\newblock 2012.
\newblock arXiv:1202.3245.

\end{thebibliography}

\end{document}